\pdfoutput=1
\RequirePackage{fix-cm}
\documentclass[12pt]{amsart}
\usepackage{mathtools}
\usepackage[all]{xy}
\usepackage{amsmath, amssymb, amsfonts, amstext, amsthm}
\usepackage[mathscr]{euscript}
\usepackage{enumerate}
\usepackage[hyphens]{url}
\newcommand{\curv}{\mathcal{C}}
\newcommand{\C}{\widetilde{\mathcal{C}}}
\newcommand{\J}{\mathcal{J}(\widetilde{\mathcal{C}})}
\newcommand{\K}{\mathcal{K}(\widetilde{\mathcal{C}})}
\newcommand{\mcO}{\mathcal{O}}
\newcommand{\F}{F_{C,A}}
\newcommand{\E}{E_{C,A}}
\newcommand{\YF}{F_{C',A'}}
\newcommand{\Hca}{H^0(C,A)}
\newcommand{\YH}{H^0(C',A')}
\newcommand{\pc}{\pi_C}

\newcommand{\Y}{\widetilde{Y}}
\newcommand{\G}{\mathcal{G}^r_d(\textrm{sm}|L|)}
\newcommand{\justg}{\mathcal{G}}
\newcommand{\mcL}{\mathcal{L}}

\newcommand{\mcV}{\mathcal{V}}

\newcommand{\wtF}{\widetilde{F}}

\newcommand{\trm}{\textrm}
\newcommand{\lr}{\longrightarrow}

\DeclarePairedDelimiter\floor{\lfloor}{\rfloor}
\newtheorem{kummer}{Theorem}[section]
\newtheorem{Theorem2}[kummer]{Theorem}
\newtheorem*{amplenef}{Remark}

\newtheorem{TheoremML}{Theorem}[section]

\newtheorem{RmkLMadapt}[TheoremML]{Remark}

\newtheorem{nefbig}[TheoremML]{Proposition}

\newtheorem{Bpfree}{Lemma}[section]
\newtheorem{BNsubset}[Bpfree]{Remark}
\newtheorem{thma}[Bpfree]{Proposition}

\newtheorem{Lisample}{Lemma}[section]

\newtheorem{smoothcurves}[Lisample]{Lemma}
\newtheorem{gonal}[Lisample]{Lemma}
\newtheorem{gononkummer}[Lisample]{Corollary}

\newtheorem{Theorem1}[Lisample]{Proposition}

\newtheorem{openset}[Lisample]{Remark}
\newtheorem{Rmkfinal}[Lisample]{Remark}
\newtheorem{g-1case}[Lisample]{Proposition}
\begin{document}
\title[\resizebox{5.5in}{!}{On the Semistability of certain Lazarsfeld-Mukai bundles on Abelian surfaces}]{On the Semistability of certain Lazarsfeld-Mukai bundles on Abelian surfaces}
\author[P. Narayanan]{Poornapushkala Narayanan}
\address{Department of Mathematics, Indian Institute of Technology Madras, Chennai - 600036.}
\email{poorna.p.narayanan@gmail.com}
\thanks{Mathematics Classification numbers: 14C20, 14C21, 14J60, 14.51}
\thanks{The final publication is available at Springer via \url{http://dx.doi.org/10.1007/s11565-017-0277-z}.}
\begin{abstract}
Let $X=\J$ be the Jacobian of a genus 2 curve $\C$ over $\mathbb{C}$ and $Y$ be the associated Kummer surface. Consider an ample line bundle $L=\mcO(m\C)$ on $X$ for an even number $m$, and its descent to $Y$, say $L'$. We show that any dominating component of $\mathcal{W}^1_{d}(|L'|)$ corresponds to $\mu_{L'}$-stable Lazarsfeld-Mukai bundles on $Y$. Further, for a smooth curve $C\in |L|$ and a base-point free $g^1_d$ on $C$, say $(A,V)$, we study the $\mu_L$-semistability of the rank-2 Lazarsfeld-Mukai bundle associated to $(C,(A,V))$ on $X$. Under certain assumptions on $C$ and the $g^1_d$, we show that the above Lazarsfeld-Mukai bundles are $\mu_L$-semistable.
\end{abstract}
\maketitle
%
\section{Introduction}\label{intro}
The kernel bundles or syzygy bundles over $\mathbb{C}$ were first studied by Ramanan and Paranjape \cite{RP}, in the context of curves. More generally, Ein and Lazarsfeld  \cite{LEM} studied these bundles over higher dimensional varieties. These bundles occur naturally as the kernels of evaluation maps of globally generated sheaves over smooth, projective varieties. In particular, let $X$ be a smooth, projective variety and $A$ be a globally generated line bundle over $X$. Then, the kernel bundle (which we denote by $M_A$) occurs as the kernel of the surjective evaluation map $H^0(X,A)\otimes\mcO_X\lr A$, i.e.,
\[ 0\lr M_A\lr H^0(X,A)\otimes\mcO_X\lr A\lr0\,. \]
The (semi)stability properties of the bundle $M_A$ (and its dual $M_A^{\vee}$) are of considerable interest. The stability of $M_A^{\vee}$ is equivalent to that of $\phi_A^*\,T_{\mathbb{P}^r}$, where $\phi_A:X\lr\mathbb{P}^r$ is the morphism associated to the globally generated line bundle $A$.

The situation has been studied extensively and is well understood in several cases, some of which are as follows:
\begin{itemize}
 \item When $X$ is a smooth curve of genus $g\geq 1$, $M_A$ is stable as soon as $\trm{deg }A\geq2g+1$, cf. \cite{RP}, \cite{BUT}, \cite{EM}, \cite{LE}, \cite{AB}, \cite{CC}.
 \item When $X=\mathbb{P}^n$ and $A=\mcO_{\mathbb{P}^n}(d)$, Flenner \cite{HF} proved the stability of $M_A$. 
 \item Camere \cite{CC1} proved the stability of kernel bundles $M_A$ when $X$ is a $K3$ surface and when $X$ is an abelian surface.
 \end{itemize}
Lazarsfeld-Mukai bundles (LM bundles) are generalizations of the (dual) kernel bundles. Suppose $X$ is a smooth projective surface over $\mathbb{C}$ and $C$ is a smooth curve on $X$. Let $A$ be a line bundle on $C$ and $V\subset H^0(C,A)$ be a linear subspace such that the linear system $\mathbb{P}V$ is base-point free. Then, $j_*A$ is a globally generated sheaf on $X$ where $j:C\hookrightarrow X$ is the inclusion. We then get a surjective evaluation map $V\otimes\mcO_X\lr j_*A$, whose kernel we denote by $F_{C,A,V}$. Thus, we get:
\[0\lr F_{C,A,V}\lr V\otimes\mcO_X\lr j_*A\lr 0\,.\]
The dual bundle $F_{C,A,V}^{\vee}=:E_{C,A,V}$ is called the \emph{Lazarsfeld-Mukai bundle}. The definition and some initial properties of LM bundles are discussed in \cite{AP} for  a $K3$ surface $X$.

Consider a $K3$ surface $X$ and an ample line bundle $L$ on $X$. Then, for a general $C\in\textrm{sm}|L|$ (cf. $\mathcal{x}\,$\ref{defns}, (b)) and a general complete, base-point free $A\in G^1_d(C)$ (cf. $\mathcal{x}\,$\ref{BN}, I.(b)), Lelli-Chiesa \cite[Chapter III]{MC} has studied the $\mu_L$-stability of the rank-2 LM bundle. Our aim in this paper is to obtain similar results on Jacobian surfaces.

In $\mathcal{x}\,$\ref{Absursect}, we study LM bundles over the Jacobian of a genus 2 curve $\C$, $X=\J$, and the associated (singular) Kummer surface  $Y=\K$. Fix an embedding of $\C\hookrightarrow X$ such that the involution $i$ on $X$ which maps $x\mapsto -x$, restricts to the hyperelliptic involution on $\C$. Now, consider the line bundle $L=\mcO(m\C)$ on $X$ for an even number $m$. Then $L$ is an ample and globally generated line bundle on $X$. Also $L$ is a symmetric line bundle, hence the involution $i$ acts on the space of global sections $H^0(X,L)$. We thus have a decomposition into the $'+'$ and $'-'$ eigenspaces: $$H^0(X,L)=H^0(X,L)^+\oplus H^0(X,L)^{-}.$$
These give $i$-fixed point linear systems $|L|^{\pm}\subset |L|$. 
Moreover, $L$ descends to a unique very ample line bundle $L'$ on the Kummer surface $Y$. 
Let  $k'$ denote the gonality (cf. $\mathcal{x}\,$\ref{cldim}, (a)) of any curve in $\textrm{sm}|L'|$ which lies in the smooth locus of $Y$. 
Then we prove the following about LM bundles on $Y=\K$.
\begin{kummer}
\label{thmkumm}
Consider a general curve $C'\in\emph{sm}|L'|$ which avoids the 16 singular points of $Y$. 
Let $d\in\mathbb{N}$ be such that $d>(\frac{m^2}{2}+1)-k'+2$ . Then, for a general $A'\in G^1_d(C')$, which is complete and base-point free, the associated LM bundle on $Y$ is $\mu_{L'}$-stable. Hence any dominating component of $\mathcal{W}^1_d(|L'|)$ (cf. $\mathcal{x}\,$\ref{BN}, III) corresponds to $\mu_{L'}$-stable LM bundles.
\end{kummer}
Note that any smooth curve $C'\in \trm{sm}|L'|$ avoiding the singular points of $Y$ has genus $g'=\frac{m^2}{2}+1$. Hence the statement requires:
$$d>g'-k'+2\,.$$
The inequality $g'-k'+2<d$ automatically implies that the Brill-Noether number (cf. $\mathcal{x}\,$\ref{BN}, II) $\rho(g',1,d)>0$, guaranteeing the non-emptiness of the Brill-Noether varieties $G^1_d(C')$ and $W^1_d(C')$ for a general $C'\in |L'|$ in the Kummer surface $Y$. In order to establish the $\mu_{L'}$-stability of LM bundles on $Y$, we look at the minimal resolution of $Y$ which is a $K3$ surface, say $\Y$, cf. \cite[Proposition VIII.11]{Be}. Let $\phi:\Y\lr Y$ be the resolution map. We then have that $H^0(\Y,\phi^*L')\simeq H^0(Y,L')$. Thus if $C'\in \textrm{sm}|L'|$ is a curve which avoids the singular points of $Y$, then $C'$ is isomorphic to a curve in $\textrm{sm}|\phi^*L'|$ in $\Y$. The result of Lelli-Chiesa which we adapted for the case of nef, big and globally generated line bundles in $\mathcal{x}\,$\ref{LMnef} tells us that for a general $C'\in \textrm{sm}|\phi^*L'|$ and a general complete and base-point free $A'\in G^1_d(C')$ in $\Y$, the LM bundle associated to $(C',A')$ is $\mu_{\phi^*L'}$-stable. Hence the corresponding bundle on $Y$ is $\mu_{L'}$-stable. Consequently, we also obtain the following on $X=\J$.
\begin{Theorem2}\label{intromain}
Let $X$, $L$ and $k'$ be as above. Let $d\in\mathbb{N}$ be such that $(\frac{m^2}{2}+1)-k'+2<d\leq \frac{m^2}{2}+2\,$.  Then there is at least one dominating component $\mathcal{W}$ of $\mathcal{G}^1_{2d}(|L|^{+})$ (cf. $\mathcal{x}\,$\ref{BN}, III) corresponding to $\mu_L$-semistable LM bundles. That is, for a general $(C,A,V)\in\mathcal{W}$, the corresponding LM bundle $F_{C,A,V}^{\vee}$ is $\mu_L$-semistable. 
\end{Theorem2}
The conditions in the statement of the theorem read as :
$$g'-k'+2<d\leq g'+1.$$
The inequality $d\leq g'+1$ tells us that for a general $C'\in |L'|$, no irreducible component of $W^1_d(C')$ is completely contained in $W^2_d(C')$ and hence in any component of $W^1_d(C')$, a general pencil is complete.

The outline of our proof is the following. We have a rational morphism 
\begin{displaymath}
\Phi:\mathcal{W}^1_d(|L'|)\dashrightarrow\mathcal{G}^1_{2d}(|L|^{+})
\end{displaymath}
where the map is defined on the open subset $\mathcal{Y}\subset\mathcal{W}^1_d(|L'|)$, given by:
 $$\mathcal{Y}=\{(C',A')\in\mathcal{W}^1_d(|L'|):C'\in\trm{sm}|L'|\trm{ avoids the 16 singular points}, h^0(C',A')=2\}.$$
The map is given by $\Phi((C',A'))=(\pi^{-1}(C')=:C,\pi|_C^*A',H^0(C',A'))$, for all $(C',A')\in\mathcal{Y}$.
Since $H^0(Y,L')\simeq H^0(X,L)^+$, the linear system of curves $|L|^{+}$ corresponds to the linear system of curves  $|L'|$. 
Hence, if $D$ is a dominating component of $\mathcal{W}^1_d(|L'|)$, then
\begin{itemize}
 \item $D$ corresponds to $\mu_{L'}$-stable LM bundles by Theorem \ref{thmkumm}, and,
 \item $\overline{\Phi(D)}$ (being irreducible and dominating) is contained in a dominating component, say $\mathcal{W}$, of $\mathcal{G}^1_{2d}(|L|^{+})$.
\end{itemize}
The quotient morphism $\pi:X=\J\lr Y=\K$ is a finite morphism of degree 2. The dual LM bundle $\YF$ on $Y$ corresponding to $(C',A')\in\mathcal{Y}$ with $A'$ base-point free pulls back to the dual LM bundle $F_{C,A,V}$ corresponding to the triple $(C, A, V)=\Phi(C',A')\in\mathcal{G}^1_{2d}(|L|^+)$ on $X$. 
Hence, by \cite[Lemma 1.17]{Mar}, the $\mu_{L'}$-semistability of LM bundles corresponding to pairs $(C',A')$ on $Y$ implies the $\mu_L$-semistability of LM bundles corresponding to triples $(C, A, V)=\Phi(C',A')$ on $X$. Thus, a general element of $\Phi(D)$ corresponds to a $\mu_L$-semistable LM bundle. Since, by \cite[Proposition 2.3.1]{HL} semistability is an open condition, this gives us a dominating component $\mathcal{W}$ of $\mathcal{G}^1_{2d}(|L|^{+})$ corresponding to $\mu_L$-semistable LM bundles. This proves the Theorem.

Finally, we consider any ample line bundle  $L$ on $X=\J$ such that the genus of a smooth curve $C\in\trm{sm}|L|$ is $g$. Then we show that a general element of $\mathcal{G}^1_{g-1}(|L|)$ corresponds to a $\mu_L$-semistable LM bundle in Proposition \ref{g-1case}. We do this by showing that the restriction of the dual LM bundle $F_{C,A,V}$ associated to a general $(C,A,V)\in\mathcal{G}^1_{g-1}(|L|)$ to $C$ i.e., $F_{C,A,V}|_C$ is $\mu_{L|_C}$-semistable.
\subsection*{Acknowledgements}
I thank Dr. Jaya NN Iyer, IMSc Chennai for introducing me to this problem and for her guidance at every stage of this project. I also thank Dr. T. E. Venkata Balaji, IIT Madras and Prof. D. S. Nagaraj, IMSc Chennai for helpful discussions and their support. I also thank Prof. Robert Lazarsfeld and Prof. Jason Starr at Stony Brook University for useful comments on an earlier version of the manuscript.

\section{Preliminaries}
\subsection{Definitions and Notations}\label{defns}
Let $X$ be a smooth projective surface over $\mathbb{C}$.  
\begin{enumerate}
 \item[(a)] If $F$ is a coherent sheaf over $X$ and $\xi$ is the generic point of $X$, then the rank of $F$ is the rank of $F_{\xi}$ as an $\mcO_{X,\xi}$ vector space.
  \item[(b)] If $L$ is a line bundle on $X$, then $\textrm{sm}|L|$ denotes the set of all curves in $|L|$ which are smooth as closed subschemes of $X$. By Bertini's theorem, we know that if $L$ is globally generated, then almost every element of $|L|$ is smooth as a closed subscheme of $X$ \cite[Corollary III.10.9]{RH}.
  \item[(c)] A line bundle $L\in \textrm{Pic}(X)$ is called nef if $(L\cdot C) \geq 0$ for all closed curves $C\subset X$. 
  \item[(d)] A line bundle $L\in \textrm{Pic}(X)$ is called big and nef if $(L\cdot L) > 0$ and L is nef. 
  By the Nakai-Moishezon Criterion \cite[Theorem V.1.10]{RH}, any ample line bundle is big and nef. 
\end{enumerate}
\subsection{(Semi)Stability with respect to nef line bundles}
\label{nefpolarization}
We recall the definition of (semi)stability of torsion-free coherent sheaves with respect to \emph{nef} line bundles. We refer to \cite{GKP} for details.

Let $L$ be a \emph{nef} line bundle on a smooth projective surface $X$ over $\mathbb{C}$. Let $E$ be a torsion-free, coherent sheaf on $X$ of rank $r>0$.
\begin{enumerate}
 \item[(a)] The slope of $E$ with respect to $L$ is defined as
 \[\mu_L(E)=\frac{c_1(E)\cdot c_1(L)}{r}\ .\]
\item[(b)] (Semi)stability with respect to a nef $L$: The torsion-free sheaf $E$ on $X$ is said to be $\mu_L$-semistable (resp.~$\mu_L$-stable) if, for any coherent subsheaf $0 \neq F \subset E$ with $\trm{rank}\,F < \trm{rank}\,E$, one has $\mu_L ( F ) \leq \mu_L ( E )$ (resp.~$\mu_L ( F ) <\mu_L ( E )$).
 \end{enumerate}
\begin{amplenef}
 Note that in the particular case when $L$ is ample, the above is just the definition of the familiar Mumford-Takemoto (semi)stability of coherent sheaves. In fact, all the basic properties of Mumford-Takemoto (semi)stability with respect to an ample line bundle, such as the existence of a unique ``Harder-Narasimhan'' filtration for torsion-free sheaves, existence of a ``Jordan-Holder'' filtration for semistable sheaves etc.  carry through naturally for the notion of (semi)stability with respect to a nef line bundle, cf. \cite[$\mathcal{x}\,2$]{GKP}. 
\end{amplenef}
\subsection{Gonality and Clifford dimension}\label{cldim}
We recall the definitions of gonality and Clifford dimension of a smooth projective curve $C$, cf. \cite{ELMS}.
\begin{itemize}
 \item[(a)] The gonality of $C$ is the smallest integer $d$ such that $C$ admits a degree $d$ (non-constant) morphism to $\mathbb{P}^1$. Note that $C$ has gonality 1 if and only if the genus of $C$ is 0.
 \item[(b)] For a line bundle $A$ on $C$, the Clifford index is defined as 
 $$\trm{Cliff }A=\trm{deg }A-2(h^0(A)-1).$$
 We say that a line bundle $A$ \emph{contributes} to the Clifford index of $C$ if $h^i(A)\geq 2$ for $i=0,1$, where the Clifford index of $C$ is
 $$\trm{Cliff }C=\trm{min}\{\trm{Cliff }A\,|\,A\in\trm{Pic}(C)\trm{ contributes to the Clifford index of }C\}.$$
 We say that $A\in \trm{Pic}\,C$ \emph{computes} $\trm{Cliff}\,C$ if $A$ contributes to $\trm{Cliff }C$ and $\trm{Cliff }A$ is minimum. The Clifford dimension of $C$ is \[\trm{min}\{h^0(A)-1\,|\,A\in\trm{Pic}(C)\,,\,A \trm{ computes the Clifford index of }C\}.\] 
\end{itemize}
\subsection{Facts from Brill-Noether theory} \label{BN}
We recall some results on Brill-Noether theory for a smooth curve $C$, cf. \cite[Chapter IV]{ACGH}. 
\begin{itemize}
 \item[(I)] 
Suppose $C$ is a smooth curve of genus $g$, we associate certain varieties to $C$:
\begin{enumerate}
 \item[(a)] The subvariety $W^r_d(C)\subset\trm{Pic}^d(C)$ parametrizes complete linear series of degree $d$ and dimension at least $r$. Set theoretically,
 $$ W^r_d(C)=\{ A\in \textrm{Pic}^d(C)\,|\,h^0(C,A)\geq {r+1}\}.$$  
 \item[(b)] The variety $G^{\,r}_d(C)$ parametrizes linear series (possibly incomplete) of degree $d$ and dimension exactly $r$ on $C$. An element of $G^{\,r}_d(C)$ is called a $g^r_d$. Set theoretically, $$G^{\,r}_d(C)=\{(A,V)\,|\,A\in \textrm{Pic}^d(C), V \subset H^0(C,A)\textrm{ is an } (r+1) \textrm{ dimensional subspace}\}. $$
A $g^r_d$, say $(A,V)$, is said to be \emph{complete} if $V=H^0(C,A)$.
\end{enumerate}
\item[(II)] If $d, r$ are integers such that $d \geq 1, r\geq 0$ and if the Brill-Noether number 
\begin{equation*}\label{nocall1}
\rho(g,r,d)=:g-(r+1)(g-d+r)\geq 0\,,
\end{equation*}
then $G^{\,r}_d(C)$ and hence $W^r_d(C)$ are non-empty. Furthermore, every component of $G^{\,r}_d(C)$ has dimension at least equal to $\rho$; the same is true for $W^r_d(C)$ provided $d\leq r+ g$.
\end{itemize}
We next recall some facts about Brill-Noether theory for a family of smooth curves, cf. \cite[Chapter 21]{ACGH2}. Let $p:\mathcal{C}\lr S$ be a family of smooth curves of genus $g>1$ parametrized by a scheme $S$. We denote $C_s=p^{-1}(s)$ for $s\in S$.
\begin{enumerate}
 \item[(III)] The following $S$-varieties can be associated to the family of smooth curves $p$.
\begin{itemize}
 \item[(a)] The variety $\mathcal{P}ic^d(p)$ parametrizes line bundles $A$ of degree $d$ on smooth curves $C_s$ for $s\in S$, and its support is given by
 $$\textrm{Supp\ } \mathcal{P}ic^d(p)=\{(s,A) \,|\, s\in S,\ A\in \textrm{Pic}^d(C_s)\}\, .$$
 \item[(b)] The subvariety $\mathcal{W}^r_d(p)\subset\mathcal{P}ic^d(p)$ whose support is given by the set 
 $$\textrm{Supp\ } \mathcal{W}^r_d(p)=\{(s,A) \,|\, s\in S,\ A\in W^r_d(C_s)\}\, ,$$
parametrizes degree $d$ line bundles $A$ on smooth curves $C_s$ for $s\in S$ with ${h^0(A)\geq r+1}$.
 \item[(c)] The variety $\mathcal{G}^r_d(p)$ whose support is given by
 $$\textrm{Supp\ }\mathcal{G}^r_d(p)=\{(s,(A,V))\,|\,s\in S,\ (A,V)\in G^r_d(C_s)\}\, ,$$
 parametrizes $g^r_d$'s on smooth curves $C_s$, for $s\in S$.
\end{itemize}
\end{enumerate}
\subsection{Lazarsfeld-Mukai bundles}\label{LMbundle}
Suppose that $X$ is a smooth projective surface and $C$ is a smooth curve on $X$. Let $A$ be a line bundle on $C$ and $V\subset H^0(C,A)$ be a linear subspace such that the linear system $\mathbb{P}V$ is base-point free. If $j:C\hookrightarrow S$ is the inclusion, then $j_*A$ is a globally generated sheaf on $X$. Therefore, we get a surjective evaluation map on $X$ given by $V\otimes\mcO_X\lr j_*A\lr 0\,.$ Let $F_{C,A,V}$ be the kernel. So, it fits in the exact sequence:
\begin{equation}\label{eqn0}
 0\lr F_{C,A,V}\lr V\otimes\mcO_X\lr j_*A\lr 0\,.
\end{equation}
Then, we have the following (cf. \cite[Proposition~5.2.2]{HL}) :
\begin{enumerate}
 \item[(a)] $F_{C,A,V}$ is a vector bundle,
 \item[(b)] $\trm{rank}(F_{C,A,V})=\textrm{\textrm{dim} } V$, $c_1(F_{C,A,V})=\mcO_X(-C)$ and $c_2(F_{C,A,V})=\trm{deg}(A)$ ,
 \item[(c)] $H^0(X,F_{C,A,V})=0$ .
\end{enumerate}
The dual of $F_{C,A,V}$, which we denote by $E_{C,A,V}$, is called the \emph{Lazarsfeld-Mukai bundle}. The bundle $E_{C,A,V}$ fits in the following exact sequence obtained by dualizing \eqref{eqn0}.
\begin{equation*}\label{nocall2}
 0\lr V^{\vee}\otimes\mcO_X\lr E_{C,A,V}\lr\mathcal{E}\!xt^1(j_*A,\mcO_X)\lr 0\,.
\end{equation*}
When $V=H^0(C,A)$ we denote $F_{C,A,V}$ by $\F$ and $E_{C,A,V}$ by $\E$.
\section{Constructing families of LM bundles on a smooth projective surface}\label{construct}
In this section, we construct families of LM bundles, parametrized by suitable varieties.
\subsection{Genericity of base-point free linear series}\label{generic} Consider a smooth projective surface $X$ over $\mathbb{C}$. Let $L$ be an ample line bundle on $X$ such that 
$\textrm{sm}|L|\neq \emptyset$. Consider 
\begin{equation*}\label{nocall3}
\mathcal{C}=\{(x,C)\,|\,C\in \trm{sm}|L|,\,x\in C\}\,.
\end{equation*}
Then $p:\mathcal{C}\lr \textrm{sm}|L|$ is a finitely presented, faithfully flat family of curves since, for any $C\in \textrm{sm}|L|$, the fiber $p^{-1}(C)=\{(x,C)\,|\,x\in C\}\cong C$. Let $d\geq 1,r\geq 0$. Given the family of curves $p$, we have the variety $\mathcal{G}^r_d(p)$ which we shall henceforth denote by $\G$. Set theoretically, 
\begin{equation*}\label{nocall4}
\G=\{(C,A,V)\,|\,C\in \textrm{sm}|L|, (A,V)\in G^r_d(C)\}\,. 
\end{equation*}
Suppose that $\G\neq \emptyset$. There is a natural projection $u:\G\lr \textrm{sm}|L|$. So, we can consider the fibered product $\justg$.
\begin{displaymath} 
\xymatrix{\justg\ \ar[r] \ar[d]_{p_1}  & \mathcal{C}\ar[d]^{p} \\
	  \G\ar[r]_u & \textrm{sm}|L|   }
\end{displaymath}
Since $p$ is flat, $p_1$ is a flat morphism too. Set theoretically,
\begin{equation*}\label{nocall5}
\justg\cong\{(x,C,A,V)\,|\,(C,A,V)\in\G, x\in C\} .
\end{equation*}
So the fiber of $p_1$ over a point $(C,A,V)\in\G$ is, in fact, $C$. Hence $p_1:\justg\rightarrow\G$ is a flat family of curves too. 
\begin{equation}\label{eqn1}
\trm{Denote }\mathcal{B}=\{(x,C,A,V)\in\justg\,|\,(A,V)\in G^r_d(C)\textrm{ is base-point free}\}\,. 
\end{equation}
\begin{Bpfree}
\label{base-pointfree}
 Let $X$ and $L$ be as above. If $\G$ (and hence $\justg$) is non-empty, then the set $\mathcal{B}\subset\justg$ described above is an open subset (maybe empty).
\end{Bpfree}\label{bpfree}
\textit{Proof.} (Due to Jason Starr) Recall that, $p:\mathcal{C}\lr \textrm{sm}|L|$ and $p_1:\justg\lr\G$ are flat families of curves. We want a family of curves which admits a section, so that it has a Poincar\'e line bundle over it. So, we once again consider the following fibered product.
\begin{displaymath} 
\xymatrix{\curv_{\justg}\ \ar[r] \ar[d]_{p_2}  & \mathcal{C}\ar[d]^{p} \\
	  \justg\ar[r]_{u\,\circ\, p_1} & \textrm{sm}|L|   }
\end{displaymath}
Set theoretically, 
\begin{equation*}\label{nocall6}
\curv_{\justg}=\{(x,(y,C,A,V))\,|\,(y,C,A,V)\in\justg,x\in C\}\,. 
\end{equation*}
So, the fiber over a point $(y,C,A,V)\in\justg$ is just $C$. Since $p$ is flat, $p_2$ is a flat family of curves. The universal property of fibered product gives us a morphism $\Delta:\justg\lr\curv_{\justg}$ such that $p_2\circ\Delta=id_{\justg}$, the identity on $\justg$. 
Hence, $\Delta:\justg\lr\curv_{\justg}$ gives a section of the family of curves $p_2:\curv_{\justg}\lr\justg$. So, there is a pair $(\mcL,\mcV)$ such that
\begin{enumerate}
 \item[(i)] $\mcL$ is the Poincar\'e line bundle on $\curv_{\justg}$ in $\mathcal{P}ic^d(\curv_{\justg})$ such that if $(x,C,A,V)\in\justg$, then $\mcL|_{p_2^{-1}(x,C,A,V)}\simeq A$.
 \item[(ii)] $\mcV$ is a locally free subsheaf of ${p_2}_*\mcL$ of rank $r+1$ such that if $t=(x,C,A,V)\in\justg$, then the subspace
 \[\mcV\otimes k(t)\hookrightarrow H^0({p_2}^{-1}(t),\mcL\otimes k(t))=H^0(C,A),\textrm{ is isomorphic to } V\,.\]
\end{enumerate}
By point (ii) above, we have the inclusion $\mcV\hookrightarrow {p_2}_*\mcL$. By adjointness, on $\curv_{\justg}$ we get
\begin{equation}\label{eqn2}
 p_2^*\mcV\lr\mcL\lr\mathcal{M}\lr 0\, ,
\end{equation}
where $\mathcal{M}$ is the cokernel. Consider $\mathcal{S}=\trm{Supp}(\mathcal{M})=\{s\in\curv_{\justg}\,|\,\mathcal{M}_s\neq 0\}$. This is a closed set. 
Since $p_2$ is a proper morphism, $p_2(\mathcal{S})\subset\justg$ is a closed set. 

Then the complement of $p_2(\mathcal{S})$ in $\justg$ parametrizes elements of $\justg$ with the corresponding $g^r_d$'s base-point free. Indeed, for $t=(x,C,A,V)\in\justg$,  restricting the morphism \eqref{eqn2} to the fiber $p_2^{-1}(t)\simeq C$, we get
$$V\otimes\mcO_C\lr A\lr \mathcal{M}|_{p_2^{-1}(t)}\,.$$
We now have this sequence of equivalent statements: $t=(x,C,A,V)\in\mathcal{B} \iff \mathcal{M}|_{p_2^{-1}(t)}=0$ $\iff$ the stalks of $\mathcal{M}|_{p_2^{-1}(t)}$ at all points of $p_2^{-1}(t)$ are zero $\iff$ $p_2^{-1}(t)\subset\mathcal{S}^c$. By a set theory argument we observe that, $p_2^{-1}(t)\subset\mathcal{S}^c\iff t\in p_2(\mathcal{S})^c$ (since $t\in p_2(\mathcal{S})^c\iff p_2^{-1}(t)\cap\mathcal{S}=\emptyset$).  Hence $p_2(\mathcal{S})^c$ is precisely our set $\mathcal{B}$, and thus $\mathcal{B}$ is open. $\qed$

\begin{BNsubset}\label{imageopen}
  Since $p_1:\justg\rightarrow \justg^r_d(\emph{sm}|L|)$ is a finitely presented, faithfully flat morphism, $p_1$ is an open map. Therefore, $$p_1(\mathcal{B})=\left\{(C,A,V)\in \justg^r_d(\emph{sm}|L|)\ \big|\ (A,V)\trm{ is base-point free}\right\}$$ is an open subset of $\justg^r_d(\emph{sm}|L|)$.
\end{BNsubset}
\subsection{Families of LM bundles over $\mathcal{B}$}\label{LMfamily}
Now, we describe the required family of LM bundles, as follows.
\begin{thma}
\label{family}
Let $X$ and $L$ be as above. Assume that 
the open subset $\mathcal{B}$ of $\justg$ given by equation \eqref{eqn1} is non-empty. Then there is a family of LM bundles of rank $r+1$ parametrized by $\mathcal{B}$. 
\end{thma}
\textit{Proof.} We keep notations as in the Lemma \ref{base-pointfree}. Recall the family of curves $p_2:\curv_{\justg}\rightarrow\justg$. Note that $\curv_{\justg}$ is a Cartier divisor in $X\times\justg$. 
\begin{displaymath} 
\xymatrix{\curv_{\justg}\ \ar@{^{(}->}[r]^j \ar[d]_{p_2}  & X\times \justg\ar[dl]^{q} & & p_2^{-1}(\mathcal{B})\, \ar@{^{(}->}[r]^j \ar[d]_{p_2}  & X\times \mathcal{B}\ar[dl]^{q} \\
	  \justg & & & \mathcal{B} }
\end{displaymath}
Consider, $p_2|_{p_2^{-1}(\mathcal{B})}:p_2^{-1}(\mathcal{B})\lr \mathcal{B}$. Then, $p_2^{-1}(\mathcal{B})$ is a Cartier divisor in $X\times \mathcal{B}$ as well. Restricting the sequence \eqref{eqn2} on $\curv_{\justg}$ to $p_2^{-1}(\mathcal{B})\subset\mathcal{S}^c\subset\curv_{\justg}$, we have the surjection
\[p_2^*\mcV\lr\mcL\lr 0\trm{ on }p_2^{-1}(\mathcal{B})\,.\]
Since $p_2=q\circ j$, we get $p_2^*=j^*\circ q^*$. So, the above sequence on $p_2^{-1}(\mathcal{B})$ becomes
\[j^*\circ q^*\mcV\lr\mcL \lr 0\,.\]
Thus by adjointness (and since $j$ is a closed immersion), we get $q^*\mcV\lr j_*\mcL\lr 0$ on $X\times \mathcal{B}$. Let $\mathcal{F}$ be the kernel. Hence, we get
\begin{equation}\label{eqn3}
 0\lr\mathcal{F}\lr q^*\mcV\lr j_*\mcL\lr 0\quad\trm{on }X\times \mathcal{B}.
\end{equation}
\begin{itemize}
 \item Since $p_2:\curv_{\justg}\lr\justg$ is flat and $\mcL$ is a line bundle on $\curv_{\justg}$, $\mcL$ is flat over $\justg$. Hence the sheaf $j_*\mcL$ on $X\times\justg$ is flat over $\justg$.
 \item Again, since the projection $X\times\justg\lr\justg$ is flat and since $q^*\mcV$ is a vector bundle on $X\times\justg$, $q^*\mcV$ is flat over $\justg$.
\end{itemize}
Hence, from \eqref{eqn3} we infer that $\mathcal{F}$ is flat over $\mathcal{B}$. If $t=(x,C,A,V)\in \mathcal{B}$, then $q^{-1}(t)\cong X$. We restrict the short exact sequence $\eqref{eqn3}$ to the fiber $q^{-1}(t)$ to get
\begin{equation*}\label{nocall7}
 0\lr\mathcal{F}|_{q^{-1}(x,C,A,V)}\lr q^*\mcV|_{q^{-1}(x,C,A,V)}\lr j_*\mcL|_{q^{-1}(x,C,A,V)}\lr 0\,.
\end{equation*}
That is,
\[0\lr\mathcal{F}|_{q^{-1}(x,C,A,V)}\lr V\otimes\mathcal{O}_X\lr j_*A\lr 0\,,\]
where, in this case, $j$ denotes the inclusion $j:C\hookrightarrow X$.

Hence, $\mathcal{F}|_{q^{-1}(x,C,A,V)}$ is the dual LM bundle $F_{C,A,V}$ associated to the triple $(C,A,V)$. Thus, we get that $\mathcal{F}$ on $X\times\mathcal{B}$ is a family of dual LM bundles on $X$ of rank $r+1$ parametrized by the open set $\mathcal{B}\subset\justg$.$\qed$
\section{Stability of rank-2 LM bundles over $K3$ surfaces with respect to nef and big line bundles}\label{LMnef}We start by recalling the theorem of Lelli-Chiesa.
\begin{TheoremML}\label{source} \cite[Theorem 3.5.3, Page 62]{MC} 
 Let $X$ be a smooth, projective $K3$ surface and $L\in \emph{Pic}(X)$ an ample line bundle such that a general curve $C\in \emph{sm}| L |$ has genus $g$, Clifford dimension 1 and maximal gonality $k = \floor*{\frac{g+3}{2}}$. If $\rho ( g, 1, d ) > 0$, then any dominating component of $\mathcal{W}^1_d(|L|)$ corresponds to $\mu_L$-stable LM bundles. That is, if $C$ is general in its linear system, then the rank-2 LM bundle associated with a general complete, base-point free $g_d^1$ on $C$ is $\mu_L$-stable.
\end{TheoremML}
Note that Theorem \ref{source} can be adapted to the case when the gonality $k'$  of curves $C\in \trm{sm}|L|$ is not necessarily maximal (i.e. $k'\leq k=\floor*{\frac{g+3}{2}}$), with an additional condition that ${d>g-k'+2}$.
%
\begin{RmkLMadapt}
 \label{rangenotempty}
We note that we can find complete $g^1_d$'s on curves $C$ satisfying $d>g-k'+2$. Indeed, to find a complete $g^1_d$ on $C\in \emph{sm}|L|$, we need $d\leq g+1$. Hence we want $d$ in the range $g-k'+2<d\leq g+1$. Since $L$ is ample, the gonality $k'$ of any curve $C\in \emph{sm}|L|$ is $>1$. Hence we can always find $d$ in this range. 
%
\end{RmkLMadapt}
We now extend Theorem \ref{source} when the ample line bundle $L$ is replaced by a nef, big and globally generated line bundle $L$. As we saw in $\mathcal{x}\,$\ref{nefpolarization}, we can study (semi)stability with respect to nef line bundles.
\begin{nefbig}
\label{nefandbig}
Let $X$ be a smooth, projective $K3$ surface and $L\in \emph{Pic}(X)$ a nef, big and globally generated line bundle such that a general curve $C\in \emph{sm}| L |$ has genus $g$, Clifford dimension 1 and gonality $k'$.  Let $d\in\mathbb{N}$ be such that $d>g-k'+2$ (hence $\rho(g,1,d)>0$). Then any dominating component of $\mathcal{W}^1_d(|L|)$ corresponds to $\mu_L$-stable LM bundles. That is, for a general $C\in \emph{sm}|L|$, the LM bundle associated to a general complete, base-point free $A\in G^1_d(C)$ is $\mu_L$-stable.
\end{nefbig}
\emph{Proof. }First we note that the definition and basic properties of LM bundles do not depend on the ampleness of $L$. 
Hence, Lelli-Chiesa's result follows by replacing the ample $L$ with a globally generated, nef and big $L$ by an application of Ramanujam-Kodaira vanishing theorem. We explain this as follows.

By \cite[Proposition 3.5.2]{MC}, we get that if $\mathcal{W}$ is any irreducible component of $\mathcal{W}^1_d(\trm{sm}|L|)$ corresponding to $\mu_L$-unstable LM bundles, then it satisfies
\[\trm{dim}\ \mathcal{W}\leq g+d-k'\,.\]
Using the bound $d>g-k'+2$, we get $-k'<d-g-2$. Hence,
\[\textrm{dim}\ \mathcal{W}<g+d+(d-g-2)=2d-2\,.\]
On the other hand any dominating component of $\mathcal{W}^1_d(\textrm{sm}|L|)$ has at least dimension equal to (dimension of base + dimension of fiber), i.e., $\textrm{dim}\, \textrm{sm}|L|+\textrm{dim}\, W^1_d(C)$, for $C\in\textrm{sm}|L|$.

Here, by Bertini's theorem, $\textrm{dim} \textrm{ sm}|L|=\textrm{dim} |L|=h^0(X,L)-1$. For any nef and big line bundle $L$, $h^i(X,L)=0$ for $i>0$ (by Ramanujam-Kodaira vanishing theorem). Therefore, $\textrm{dim} \textrm{ sm}|L|=\chi(L)-1=\frac{L^2}{2}+2-1=g$. Hence any dominating component of $\mathcal{W}^1_d(\textrm{sm}|L|)$ has
\[\textrm{dimension}\geq g+\rho(g,1,d)=2d-2.\]
Therefore $\mathcal{W}$ is not a dominating component of $\mathcal{W}^1_d(\textrm{sm}|L|)$. Hence any dominating component of $\mathcal{W}^1_d(\textrm{sm}|L|)$ corresponds to $\mu_L$-stable LM bundles. $\qed$
\section{LM bundles over abelian surfaces}\label{Absursect}
Consider a smooth projective curve $\C$ of genus 2 over $\mathbb{C}$. Suppose $X:=\J$, the Jacobian of $\C$. Let $i:X\lr X$ be the involution which maps $x\mapsto -x$. Consider the quotient of $X$ by the action of $i$, $Y:={X}/{(i)}$. Then $Y$ is the Kummer surface $\K$ associated to the abelian surface $X$. This is a singular surface whose singularities are 16 nodes, which are all the images of the 2-torsion points of $X$. Let $\pi:X\lr Y$ be the canonical quotient morphism of degree 2. In fact, if $U$ denotes the smooth locus of $Y$, $\pi|_{\pi^{-1}(U)}:\pi^{-1}(U)\lr U$ is a flat 2-sheeted covering morphism.

We fix an embedding $\C\hookrightarrow X$, via a Weierstrass point, such that the involution $i$ on $X$ restricts to the hyperelliptic involution on $\C$.  
Therefore, we have $i(\C)=\C$ and $i^*\mcO(\C)=\mcO(\C)$ on $X$. Hence $\mcO(\C)$ is a \emph{symmetric} line bundle on $X$, cf. \cite[$\mathcal{x}$2]{Mum}. If we choose $L$ to be the line bundle $L=\mcO(m\C)$, where $m$ is an even natural number, then $L$ is a totally symmetric line bundle on the abelian surface $X$. Thus $L$ descends to the Kummer surface $Y$. For details cf. \cite[$\mathcal{x}2$, Proposition 1]{Mum}. Let $L'$ be the unique line bundle on $Y$ such that $\pi^*L'\simeq L$. 
\begin{Lisample}\label{veryample}
 Both $L$ and $L'$ are ample line bundles on $X$ and $Y$ respectively, with self-intersection numbers $(L^2)=2m^2$ and $(L'^{\,2})=m^2$. In fact, $L'$ is very ample.
\end{Lisample}
\textit{Proof.} This is an application of the Riemann-Roch theorem. Indeed, since $L=\mcO(m\C)$, $L$ is a multiple of the theta divisor and hence is ample. 
Thus, by Riemann-Roch formula, $h^0(X,L)=m^2$ and $(L^2)=m^2\C^2=2m^2$.

The descent of $\mcO(2\C)$ defines an embedding of $Y$ in $\mathbb{P}^3$, and hence is very ample \cite[Theorem IV.8.1]{BL}. Similarly, we see that $L'$ is very ample on $Y$ when $m$ is even. Also since $\pi^*L'\simeq L$, $(L^2)=2m^2$ and $\pi$ is a degree 2 morphism, we get that $(L'^{\,2})=m^2$. $\qed$

Next consider the minimal resolution of the singular surface $Y$ to get the smooth surface $\Y$. Let $\phi:\Y\lr Y$ denote the resolution map. In fact, $\Y$ is a $K3$ surface, cf. \cite[Proposition VIII.11]{Be}. Consider the line bundle $\phi^*L'$ on $\Y$. Since $\phi$ is a birational morphism, $(\phi^*L'^{\,2})=(L'^{\,2})=m^2$. But $\phi^*L'$ is \emph{not ample} because it is trivial when restricted to any exceptional divisor. However, we have the following lemma. 
\begin{smoothcurves}
\label{smcurves}
On the $K3$ surface $\Y$ we have the following.
\begin{itemize}
 \item[(a)] The line bundle $\phi^*L'$ is nef, big and globally generated.
 \item[(b)] $\emph{dim}\ |\phi^*L'| =\frac{m^2}{2}+1$.
 \item[(c)] A general member of the linear system $|\phi^*L'|$ is smooth.
 \end{itemize}
 \end{smoothcurves}
\textit{Proof.} By Lemma \ref{veryample}, $L'$ is very ample, therefore it is globally generated, nef and big. Hence, the pull back $\phi^*L'$ is also globally generated, nef and big. This proves (a). Since $\phi^*L'$ is nef and big, by Ramanujam-Kodaira vanishing theorem we get that $h^i(\Y,\phi^*L')=0$ for $i=1,2$. Thus, by Riemann-Roch theorem, we get \[h^0(\Y,\phi^*L')=\chi(\phi^*L')=\frac{(\phi^*L'^{\,2})}{2}+2=\frac{m^2}{2}+2\,.\] Hence, $\textrm{dim}\ |\phi^*L'|=\frac{m^2}{2}+1$. Since $\phi^*L'$ is a globally generated line bundle on $\Y$, by Bertini's theorem, a general member $C\in |\phi^*L'|$ is smooth. $\qed$

 We also make the following observations about the nef, big and globally generated line bundle $\phi^*L'$ on $\Y$.
\begin{gonal}\label{gonal}
The smooth curves in the linear system $|\phi^*L'|$ have genus $\frac{m^2}{2}+1$, Clifford dimension 1 and the same gonality, say $k'$.
\end{gonal}
\textit{Proof.} Any curve in $\textrm{sm}|\phi^*L'|$ has genus $g'=\frac{m^2}{2}+1$. For, by adjunction formula, $m^2=(\phi^*L'^{\,2})=2g'-2$.

The gonality of any smooth curve in the linear system is constant, say $k'$. This is because Knutsen \cite[Theorem 1]{Kn} proved that if $S$ is a $K3$ surface and $E$ is a globally generated line bundle on S, and if the gonality of curves in $\trm{sm}|E|$ is not constant, then the pair $(S,E)$ is as in Donagi and Morrison's example \cite{DM}. In their example, $f:S\lr\mathbb{P}^2$ is a $K3$ surface which is a double cover of $\mathbb{P}^2$ branched along a nonsingular plane sextic curve $B\subset\mathbb{P}^2$ and $E$ is the line bundle $f^*\mathcal{O}_{\mathbb{P}^2}(3)$. But our $K3$ surface $\Y$ and our globally generated line bundle $\phi^*L'$ are not of this form. Hence we get the constancy of gonality.

Again by \cite{Kn}, except for the Donagi-Morrison example, if a globally generated linear system $|E|$ on a $K3$ surface $S$ contains smooth curves of Clifford dimension at least 2, then $E = \mcO_S(2D+\Gamma)$, where $D, \Gamma \subset S$  are smooth curves, $D^2 \geq 2$, $\Gamma^2 = -2$, and $D\cdot\Gamma=1$, cf. \cite[ $\mathcal{x}\,5$]{AF}. We have already eliminated the possibility that $(\Y,\phi^*L')$ is the Donagi-Morrison example. The other case is that $\phi^*L'$ has the form $\phi^*L'=\mcO_{\Y}(2D+\Gamma)$. This implies that 
\[(2D+\Gamma)^2=(\phi^*L'^{\,2})=m^2\textrm{ , i.e.\,,}\]
\[4D^2+\Gamma^2+4D\cdot\Gamma=m^2,\textrm{ and hence\,,}\]
\[4D^2+2=m^2\,.\] 
This is again not possible since $m$ is even. Hence we have eliminated the second possibility too. Therefore curves in $\textrm{sm}|\phi^*L'|$ have Clifford dimension 1.$\qed$

Since the $K3$-surface $\Y$ is isomorphic to the Kummer surface outside the exceptional loci, we have the following easy consequence.
\begin{gononkummer}
 \label{gonalonkumm}
Consider the line bundle $L'$ on $Y=\K$, as in Lemma \ref{veryample}. Then, any curve $C'\in \emph{sm}|L'|$ which avoids the 16 singular points of $Y$ has genus $\frac{m^2}{2}+1$, (constant) gonality $k'$ (from Lemma \ref{gonal}) and Clifford dimension 1. 
\end{gononkummer}
\begin{openset}
 Note that, the curves in $\emph{sm}|L'|$ on $Y$ which avoid the 16 singular points is an open subset of $|L'|$. So also, the curves in $\emph{sm}|\phi^*L'|$ on $\Y$ which avoid the exceptional loci is an open subset of $|\phi^*L'|$. 
\end{openset}
\subsection{LM bundles on $Y$ and $\Y$}\label{LMkummer}
We begin with  a proof of Theorem \ref{thmkumm}. We now see that $\Y$ is a $K3$ surface and $\phi^*L'$ is a nef, big and globally generated line bundle on $\Y$ such that a general curve in $\textrm{sm}|\phi^*L'|$ has genus $g'$ (=$\,\frac{m^2}{2}+1$), gonality $k'$ and Clifford dimension 1. Let $d\in\mathbb{N}$ be such that $d>g'-k'+2$ (thus $\rho(g',1,d)>0$). Then by Proposition \ref{nefandbig}, any dominating component of $\mathcal{W}^1_d(|\phi^*L'|)$ corresponds to $\mu_{\phi^*L'}$-stable LM bundles. 

\textbf{Proof of Theorem \ref{thmkumm}.} Fix any $d\in\mathbb{N}$ satisfying the given conditions. 
Now, consider a curve $C'\in \textrm{sm}|L'|$ which avoids the 16 singular points of $Y$. 
Let $A'$ be a complete, base-point free $g^1_d$  on $C'$. Since $\phi:\Y\lr Y$ is an isomorphism outside the exceptional loci, we have the following commutative diagram.
\begin{displaymath}
\xymatrix{A'\ar[r]\ar[d] & C'\ \ar@{^{(}->}[r]^{j'} \ar[d]^{\phi|_{C'}}_{\sim} & \Y\ar[d]^\phi\\
	  A'\ar[r] & C'\ar@{^{(}->}[r]_{j'} & Y }
\end{displaymath}
We have the following two exact sequences with respect to the dual LM bundles (which are of rank-2) corresponding to pairs $(C',A')$ on $Y$ and $(C',A')$ on $\Y$, say $F$ and $\widetilde{F}$ respectively:
\begin{equation}\label{eqn4}
 0\lr F\lr \YH\otimes\mcO_{Y}\lr j'_*A'\lr 0\, ,
\end{equation}
\begin{equation}\label{eqn5}
 0\lr \widetilde{F}\lr \YH\otimes\mcO_{\Y}\lr j'_*A'\lr 0\,.
\end{equation}
Clearly the exact sequence \eqref{eqn5} is the pull back of the exact sequence \eqref{eqn4} under $\phi$. Hence $\phi^*F=\widetilde{F}$.

As remarked earlier, for a general $C'\in \textrm{sm}|\phi^*L'|$ and a general complete, base-point free $A'\in G^1_d(C')$ the corresponding dual LM-bundle $\widetilde{F}$ is $\mu_{\phi^*L'}$-stable. We then claim that $F$ is $\mu_{L'}$-stable.

Suppose $G\subset F$ is a line subbundle of $F$. Note that $\phi^*G\subset\phi^*F=\wtF$ is again a line subbundle of $\wtF$. For, $\phi^*G$ is a line bundle, and a morphism from a line bundle to any vector bundle is either injective or zero. But we know that this morphism is injective on an open set (complement of the exceptional loci), and hence is non-zero.

Now since $\wtF$ is $\mu_{\phi^*L'}$-stable,
\begin{equation*}\label{nocall8}
c_1(\phi^*G)\cdot(\phi^*L')<\frac{c_1(\phi^*F)\cdot(\phi^*L')}{2}\,.
\end{equation*}
Since $\phi$ is the blow up map, $c_1(\phi^*F)\cdot(\phi^*L')=c_1(F)\cdot L'$, and similarly $c_1(\phi^*G)\cdot(\phi^*L')=c_1(G)\cdot L'$. Hence, the above inequality becomes,
\[c_1(G)\cdot L'<\frac{c_1(F)\cdot L'}{2}\,.\]
This shows that $F$ is $\mu_{L'}$-stable. This proves our theorem. $\qed$

Now we resume our analysis on $X=\J$.

 As earlier, choose $m$ to be even. Then, $L=\mcO(m\C)$ is an ample and globally generated line bundle on $X$ (in fact $L$ is very ample for $m\geq 4$, cf. \cite[Corollary IV.5.3]{BL}). Therefore, by Bertini's theorem $\trm{sm}|L|\neq\emptyset$. Recall that $i^*L=L$. Hence the involution $i$ acts on $H^0(X,L)$, and we have a decomposition into the $'+'$ and $'-'$ eigenspaces: $$H^0(X,L)=H^0(X,L)^+\oplus H^0(X,L)^{-}.$$
These give $i$-fixed point linear systems $|L|^{\pm}\subset |L|$. We have dim$|L|^{\pm}=h^0(X,L)^{\pm}-1$, where by \cite[Corollary IV.6.6]{BL}, 
\begin{equation*}\label{nocall9}
h^0(X,L)^{\pm}=\frac{1}{2}h^0(X,L)\pm 2=\frac{m^2}{2}\pm 2\,.
\end{equation*}
Let $C\in\textrm{sm}|L|$ be a smooth curve in $X$, and $C'=\pi(C)$ be the image in the Kummer surface $Y$. Denote $\pi|_C=:\pi_C$.
\begin{enumerate}
 \item[(a)] If $C\notin |L|^{\pm}$, then $C$ is not preserved under the involution. Yet, we have that $i(C)\cap C\neq\emptyset$. For, since $C$ is an ample divisor, by the Nakai-Moishezon Criterion, $C\cdot i(C)>0$. In such a case we have finitely many points $x\in C$ such that $-x\in C$ as well. Then the image $C'$ in the Kummer surface $Y$ is singular.
 \item[(b)] The other case is that $C$ is preserved under the involution. Hence $C\in |L|^{\pm}$. Then $C'$ is the quotient of the smooth curve $C$ by the action of the finite group $\{1,i\}$. So $C'$ is smooth. 
 In this case $\pc:C\lr C'$ is a finite double cover. Further, if $C\in |L|^{-}$, it passes through the 16 double points of $X$ since these are the base-points of the linear system $|L|^{-}$.
\end{enumerate}
Note that, since $L$ is a totally symmetric line bundle, 
\[H^0(X,L)^+=H^0(Y,L')\,.\]
So curves in $|L'|$ correspond to curves in $|L|^+$. Hence we restrict our attention to the linear system $|L|^+$, noting that a $C\in \textrm{sm}|L|^+$ is a smooth curve preserved under the involution and a general such $C$ avoids the 16 double points of $X$. 

Now, choose $C'\in \textrm{sm}|L'|$ which avoids the 16 singular points in $Y$. Let $C=:\pi^{-1}(C')$. 
Hence $C\in \trm{sm}|L|^+$. Let $A'$ be a globally generated line bundle on $C'$. This pulls back to a globally generated line bundle $A$ on $C$, that is, $A=\pi_C^*\,A'$.   
\begin{displaymath}
\xymatrix{A\ar[r]\ar[d] & C\ \ar@{^{(}->}[r]^j \ar[d]^{\pi_C} & X\ar[d]^\pi\\
	  A'\ar[r] & C'\ar@{^{(}->}[r]_{j'} & Y }
\end{displaymath}
Since $\pi_C:C\lr C'$ is a degree 2 \'etale morphism of smooth curves, $$\trm{deg}(A)=\trm{deg }({\pi_C}^*A')=\trm{deg}(\pi_C)\trm{deg}(A')=2\cdot\trm{deg}(A')\,.$$ 
Let $g(C)$ denote the genus of $C$. By adjunction formula, any curve in $ \textrm{sm}|L|$ has genus ${\frac{(L^2)}{2}+1}$. Now, $L=\mcO(m\C)$, therefore $(L^2)=2m^2$ and $g(C)=m^2+1$. Also, by Riemann-Hurwitz formula $2g(C)-2=2(2g(C')-2)$. So, $g(C')=\frac{m^2}{2}+1$. Note that ${\pc}_{*}\mathcal{O}_C$ is a rank-2 vector bundle on $C'$ of the form ${\pc}_{*} \mcO_C\simeq\mcO_{C'}\oplus M$, where $M\in \textrm{Pic}^0(C')$ such that $M^2=\mcO_{C'}$. Hence for $A=\pc^*\,A'$, we see by projection formula that:
\begin{eqnarray*}\label{nocall10}
{\pi_C}_*A & \simeq & {\pi_C}_*{\pi_C}^*A'  \simeq \ {\pi_C}_*\mcO_C\otimes A' {}
\nonumber\\
& \simeq & (\mcO_{C'}\oplus M)\otimes A' {}
\nonumber\\
& \simeq & A'\oplus (A'\otimes M)\,.
\end{eqnarray*}
Thereby, we see that the global sections of $A$ on $C$ can be written as:
\begin{eqnarray}\label{nocall11}
\Hca & = & H^0(C',{\pi_C}_*A)\ {}
\nonumber\\
& = & \pi_C^*H^0(C',A'\oplus A'\otimes M) {}
\nonumber\\
& = & \pi_C^*\YH\oplus \pi_C^*H^0(C',A'\otimes M)\,.
\end{eqnarray}
Let ${V=H^0(C',A')}$. Since $A'$ is globally generated on $C'$, we have ${V\otimes \mcO_{C'}\lr A'\lr 0}$. Hence pulling back by $\pi_C$, we get the surjection $V\otimes\mcO_C\lr A\lr 0$. Hence, $V=H^0(C',A')\subset H^0(C,A)$ is a subspace such that the linear system $\mathbb{P}V$ on $C$ is base-point free. Now, we consider the dual LM bundle on $X$ corresponding to $(C,A,V)$, which is given by the following exact sequence :
\begin{equation*}\label{eqn6}
 0\lr F_{C,A,V}\lr V\otimes\mcO_X\lr j_*A\lr 0\,.
\end{equation*}
We then have the following Proposition.
\begin{Theorem1}\label{prop}
 The dual bundle $F_{C,A,V}\simeq\pi^*\YF$. Moreover, $\YF$ is $\mu_{L'}$-semistable if and only if $F_{C,A,V}$ is $\mu_L$-semistable.
\end{Theorem1}
\emph{Proof:} The dual LM bundle on $Y$ associated to the pair $(C',A')$ is given by:
\begin{equation*}\label{nocall12}
 0\lr\YF\lr V\otimes\mcO_Y\lr j'_*A'\lr 0\,.
\end{equation*}
We pullback this short exact sequence by $\pi$, to get the short exact sequence on $X$: 
\begin{equation}\label{eqn7}
 0\lr\pi^*\YF\lr V\otimes\mcO_X\lr \pi^*j'_*A'\lr 0\,.
\end{equation}
Now if $U\subset Y$ denotes the smooth locus, then $\pi|_{\pi^{-1}(U)}:\pi^{-1}(U)\lr U$ is a flat morphism. We also have $j:C\hookrightarrow \pi^{-1}(U)$ and $j':C'\hookrightarrow U$. Hence, we have the following commutative diagram : 
\begin{displaymath}
\xymatrix{C\ar[r]^\pc\ar[d]^j & C'\ar[d]^{j'}\\
	  \pi^{-1}(U)\ar[r]^{\pi}_{flat} & U}
\end{displaymath}
Note that since $j'$ is a closed immersion, it is a separated and finite type morphism of noetherian schemes. Then, by \cite[Proposition III.9.3]{RH}, we have $\pi^*j'_*A'\simeq j_*\pi_C^*A'\simeq j_*A$ on $\pi^{-1}(U)$, and hence on $X$. Thus, from \eqref{eqn7}, we get the following exact sequence on $X$:
\begin{equation*}\label{nocall13}
 0\lr\pi^*\YF\lr V\otimes\mcO_X\lr j_*A\lr 0\,.
\end{equation*}
Therefore, $\pi^*\YF\simeq F_{C,A,V}$. By \cite[Lemma 1.17]{Mar}, if $f:X\lr Y$ is a finite morphism of normal projective varieties of dimension $n$ over an algebraically closed field $k$, and if $H$ is an ample line bundle on $Y$, then a vector bundle $F$ on $Y$ is $\mu_{H}$-semistable if and only if $f^*F$ is $\mu_{f^*H}$-semistable. Hence $\YF$ is $\mu_{L'}$-semistable if and only if $F_{C,A,V}$ is {$\mu_{L}$-semistable. $\qed$
\subsection{Semistability LM bundles on $X=\J$} We now go on to prove Theorem \ref{intromain}.

\textbf{Proof of Theorem \ref{intromain}}:  Recall that $k'$ denotes the gonality of any curve in $\trm{sm}|L'|$ which lies in the smooth locus of $Y$. The genus of any such curve is $g'=\frac{m^2}{2}+1>0$, by Corollary \ref{gonalonkumm}. Hence, $k'>1$. Therefore, the given range $g'-k'+2<d\leq g'+1$ is non-empty. 
 Fix a $d$ in the above range. Note that such a $d$ automatically satisfies $\rho(g',1,d)>0$ as well. 
 
 Consider the subset $\mathcal{Y}\subset\mathcal{W}^1_d(|L'|)$ given by:
 $$\mathcal{Y}=\{(C',A')\in\mathcal{W}^1_d(|L'|) :C'\in\trm{sm}|L'|\trm{ avoids the 16 singular points}, h^0(C',A')=2\}.$$
 Note that $\mathcal{Y}$ is a non-empty open subset of $\mathcal{W}^1_d(|L'|)$. For, the set $\mathcal{Y}$ 
parametrizes smooth curves in $|L'|$ lying in the smooth locus of $Y$, together with \emph{complete} $g^1_d$'s on the curves. Since, by hypothesis, $d\leq g'+1$, no component of $\mathcal{W}^1_d(|L'|)$ is completely contained in $\mathcal{W}^2_d(|L'|)$. Thus $\mathcal{W}^1_d(|L'|)\setminus\mathcal{W}^2_d(|L'|)$ and thereby $\mathcal{Y}$ is a non-empty open set.
 
 We now have a rational morphism $\Phi$ defined on $\mathcal{Y}$:
 $$\Phi:\mathcal{W}^1_d(|L'|)\dashrightarrow\mathcal{G}^1_{2d}(|L|^+)\quad :\quad (C',A')\mapsto(C,A, H^0(C',A')),$$
 where $C=\pi^{-1}(C')$, $A\simeq\pi|_C^*(A')$. By \eqref{nocall11}, $H^0(C',A')\subset H^0(C,A)$ is a subspace. Also $C\in\trm{sm}|L|^+$ because, (since $H^0(Y,L')\simeq H^0(X,L)^+$) curves in $\trm{sm}|L|^+$ precisely consists of preimages of curves in $\trm{sm}|L'|$. Hence, if $D$ is a dominating component  of $\mathcal{W}^1_d(|L'|)$, then $\overline{\Phi(D)}$ (being irreducible and dominating) is contained in a dominating component, say $\mathcal{W}$ of $\mathcal{G}^1_{2d}(|L|^+)$.
 
 
By Theorem \ref{thmkumm}, any dominating component of $\mathcal{W}^1_d(|L'|)$ corresponds to $\mu_{L'}$-stable LM bundles. If we start with a general $(C',A')$ in a dominating component, say $D$, of $\mathcal{W}^1_d(|L'|)$, where $C'\in \textrm{sm}|L'|$ avoids the 16 singular points of $Y$, and $A'$ is a \emph{complete}, base-point free $g^1_d$ on $C'$, 
then the LM bundle on $Y$ associated to $(C',A')$ is $\mu_{L'}$-stable. If $(C,A,V)=\Phi(C',A')$, then 
by Proposition \ref{prop}, the dual LM bundle $F_{C,A,V}$ on $X$ associated to the triple $(C,A,V)$ is $\mu_L$-semistable. This proves that the LM bundle on $X$ associated to a general $(C,A,V)\in\Phi(D)$ is $\mu_L$-semistable. We now claim that $\mathcal{W}$ corresponds to $\mu_L$-semistable LM bundles, where $\mathcal{W}$ is the dominating component of $\mathcal{G}^1_{2d}(|L|^+)$ which contains $\overline{\Phi(D)}$.

Consider the open subset $\mathcal{B}\subset\justg=\mathcal{G}^1_{2d}(\textrm{sm}|L|^+)\times_{\trm{sm}|L|^+}\mathcal{C}$ defined by equation \eqref{eqn1}. The open set $\mathcal{B}$ is non-empty. This is because a general element of any dominating component of $\mathcal{W}^1_d(|\phi^*L'|)$ on the $K3$ surface $\widetilde{Y}$ is of the form $(C',A')$ where $C'\in\trm{sm}|\phi^*L'|$ and $A'$ is complete and base-point free since $d>g'-k'+2$, cf. for instance Remark 3.13 in \cite{AF}. Thus the same is true for any dominating component of $\mathcal{W}^1_d(|L'|)$ in $Y=\K$. Consider a general element $(C',A')$ of $\mathcal{W}^1_d(|L'|)$ with $C'\in\trm{sm}|L'|$ lying in the smooth locus of $Y$, and $A'\in W^1_d(C')$ \emph{complete} and base-point free. Then $\Phi(C',A')=(C,A,H^0(C',A'))\in\mathcal{G}^1_{2d}(|L|^+)$ is such that $C\in\trm{sm}|L|^+$, and $(A,H^0(C',A'))\in G^1_{2d}(C)$ is base-point free. Then $(x,C,A,H^0(C',A'))\in\mathcal{B}$ for any $x\in C$ and thus $\mathcal{B}$ is a non-empty open set. 

Now, consider the family of LM bundles on $X$ parametrized by $\mathcal{B}$ that we constructed in Proposition \ref{family}. This is a family of LM bundles associated to tuples $(x,C,A,V)$ where $C\in \textrm{sm}|L|^+$, $x\in C$ and $(A,V)\in G^1_{2d}(C)$ is base-point free. In order to say that the LM bundle corresponding to a general $(C,A,V)\in \mathcal{W}$ is semistable, it is enough to produce a single $(x,C,A,V)\in\mathcal{B}$  with $(C,A,V)\in \mathcal{W}$ such that the corresponding dual LM bundle $F_{C,A,V}$ is $\mu_L$-semistable.  This is because, by \cite[Proposition 2.3.1]{HL}, semistability of coherent sheaves is an open property in flat families. If $(C,A,V)$ is a general element of $\Phi(D)\subset\mathcal{W}$, then $F_{C,A,V}$ is $\mu_L$-semistable. Thus for such a $(C,A,V)\in\mathcal{W}$, a corresponding $(x,C,A,V)$ with $x\in C$ is a candidate. This proves our claim and hence the theorem. $\qed$
\begin{Rmkfinal}\label{Rmkfinal}
In the case when $m=2$, that is when we are considering the line bundle $L=\mcO(2\C)$ on $X=\J$, note that $$H^0(X,L)=H^0(X,L)^+\,,\trm{ therefore } |L|=|L|^+\,.$$
Hence, by our Theorem \ref{intromain}, we get at least one irreducible component of $\mathcal{G}^1_{2d}(|L|)$ which dominates the linear system $|L|$ corresponding to $\mu_L$-semistable LM bundles.
\end{Rmkfinal}
We next show by other means that, if $L$ is any ample line bundle on $X=\J$, such that any curve in $\trm{sm}|L|$ has genus $g$, then the LM bundle corresponding to a general element of $\mathcal{G}^1_{g-1}(|L|)$ is $\mu_L$-semistable.
\begin{g-1case}\label{g-1case}
 On $X=\J$, let $L$ be an ample line bundle such that a general curve in $C\in\emph{sm}|L|$ has genus $g$. Then for a base-point free $(A,V)\in G^1_{g-1}(C)$, the associated dual LM bundle $F_{C,A,V}^{\vee}$ is $\mu_L$-semistable.
\end{g-1case}
\emph{Proof.} Let $C\in\trm{sm}|L|$. If $(A,V)\in G^1_{g-1}(C)$ is base-point free, we have the surjective map of vector bundles $V\otimes\mcO_C\twoheadrightarrow A$, whose kernel is the line bundle $A^\vee$. Hence we have the following short exact sequence on $C$:
\begin{equation}\label{SESforCAVonC}
0\lr A^{\vee}\lr V\otimes \mcO_C\lr A\lr 0\,.
\end{equation}
Also recall that the dual LM bundle on $X$ is given by the exact sequence \eqref{eqn0}. Restricting this short exact sequence \eqref{eqn0} from $X$ to $C$, we get the exact sequence, where $K$ is the kernel:
\begin{equation}\label{SESresttoC}
 0\lr K\lr F_{C,A,V}|_C\lr V\otimes\mcO_C\lr A\lr 0\,.
 \end{equation}
Here $K$ is a rank one torsion-free sheaf on the curve and hence a line bundle. Since determinant is additive in exact sequences, we get:
\[K\otimes \trm{det}(F_{C,A,V}|_C)^{\vee}\otimes \mcO_C\otimes A^{\vee}\simeq \mcO_C\,,\trm{ thus,}\]
\[K\simeq \trm{det}(F_{C,A,V}|_C)\otimes A\,.\]
From $\mathcal{x}\,$\ref{LMbundle}, we know that $\trm{det}\,F_{C,A,V}\simeq\mcO_X(-C)\simeq L^{\vee}$. Thus $\trm{det}\,(F_{C,A,V}|_C)\simeq L^{\vee}|_C\simeq \omega_C^{\vee}$, by adjunction formula, since $\omega_X\simeq\mcO_X$. Thus $K\simeq \omega_C^{\vee}\otimes A$. Further from the exact sequences \eqref{SESforCAVonC} and \eqref{SESresttoC}, we get the short exact sequence on $C$:
\begin{equation*}\label{}
 0\lr K\simeq\omega_C^{\vee}\otimes A\lr F_{C,A,V}|_C\lr A^{\vee}\lr 0\,.
\end{equation*}
Now, $$\mu_{L|_C}(K)=\trm{deg}\,(\omega_C^{\vee}\otimes A)=-(2g-2)+(g-1)=-(g-1)=\trm{deg}\,A^{\vee}=\mu_{L|_C}(A^{\vee}).$$
Hence $F_{C,A,V}|_C$ is a vector bundle on $C$ which is an extension of two line bundles with the same slope. Thus, by \cite[Lemma 1.10(3)]{Mar}, $F_{C,A,V}|_C$ is $\mu_{L|_C}$-semistable on $C$. This implies that $F_{C,A,V}$ is itself $\mu_L$-semistable, cf. for example \cite[Chapter 11]{LP}. $\qed$


\end{document}